\documentclass[a4paper,12pt]{article}
\usepackage[cp1251]{inputenc}
\usepackage[russian]{babel}
\usepackage{amsfonts, amssymb, amsmath, amsthm, amscd}
\usepackage{cite}
\author{A.R. Alimov}
\title{Monotone path-\allowbreak con\-nec\-ted\-ness and solarity of Menger-connected sets}
\date{}
\begin{document}

\maketitle

\newenvironment{Biblio}{%
                  \renewcommand{\refname}{\footnotesize REFERENCES}%
                  \begin{thebibliography}{99}%
                  \itemsep -0.2em}%
                  {\end{thebibliography}}

\def\inff{\mathop{\smash\inf\vphantom\sup}}
\renewcommand{\le}{\leqslant}
\renewcommand{\ge}{\geqslant}
\theoremstyle{plain}
\newtheorem{theorem}{Theorem}
\newtheorem{definition}{Definition}
\newtheorem{corollary}{Corollary}
\newtheorem{lemma}{Lemma}
\newtheorem{remark}{Remark}
\newtheorem{propos}{Proposition}
\renewcommand{\proofname}{\bf Proof}
\newtheorem{theorema}{Theorem}
\renewcommand{\thetheorema}{\Alph{theorema}}
\newtheorem{problem}{Open problem}\def\theproblem{}

\numberwithin{equation}{section}
\let\emptyset\varnothing
\def\mcc{\operatorname{m}}
\def\stone{\textit{MeI}}
\def\sttwo{\textit{Ex-}w^*\!\textit{s}}
\renewcommand{\abstractname}{Abstract}

\def\brL{[\kern-1.15pt[}
\def\brR{]\kern-1.15pt]}

\maketitle
\begin{abstract}
A~boundedly compact (boundedly weakly compact) $\mcc$-con\-nec\-ted  (Men\-ger-connected) set is shown to be
monotone path-\allowbreak con\-nec\-ted and is a~sun
in a~broad class of Banach spaces (in particular, in separable spaces).
Further, the intersection of a~boundedly compact monotone path-\allowbreak con\-nec\-ted ($\mcc$-connected) set with a~closed ball is established to be cell-like (trivial
shape) and, in particular, acyclic
(contractible, in the finite-dimensional case) and a~sun. Further, a~boundedly weakly compact $\mcc$-con\-nec\-ted  set is asserted to be monotone path-\allowbreak con\-nec\-ted.
In passing, Rainwater--\allowbreak Simons's theorem on weak convergence of sequences
is extended to the convergence with respect to the associated norm (in the sense of Brown).

Bibliography 39 titles.
\end{abstract}

\textbf{Keywords}:
sun, acyclic set, cell-like set, monotone path-\allowbreak con\-nec\-ted set, Menger-connected set, $d$-convexity, Menger-convexity, Rainwater--\allowbreak Simons convergence theorem

\markright{Solarity of  monotone path-\allowbreak con\-nec\-ted sets}

\footnotetext[0]{This research was carried out with the support of the Russian Foundation for Basic Research (grant no.~13-01-00022).}

\section{Introduction and main definitions}

For a bounded subset $M\ne \emptyset$ of a~normed linear space $X$, we let $\mcc (M)$  denote the intersection of all
closed balls containing~$M$ (following Brown~\cite{Brown1987},
 $\mcc (M)$ is referred to as the {\it Banach--\allowbreak Mazur hull} (the ball hull)
of a~set~$M$. A~set $M\subset X$ is called $\mcc$-{\it connected\/}
(\textit{Menger-connected}) \cite{Brown1987} if
$  \mcc(\{x,y\})\cap M\ne \{x,y\} $ for any $x,y\in M$. For brevity, we write $ \mcc(\{x,y\})=\mcc(x,y)$.
Despite the name, a~closed $\mcc$-connected subset of an in the infinite-dimensional need not be connected.

Let $k(\tau)$, $0\le \tau\le 1$, be a~continuous curve in a~normed linear space~$X$. Following~\cite{BDL},
we say that the curve  $k(\cdot)$ is {\it monotone} if $f(k(\tau))$ is a~monotone function of~$\tau$ for any $f\in \operatorname{ext} S^*$
(here and henceforth,  $\operatorname{ext} S^*$ is the set of extreme points of the unit sphere  $S^*$ of the dual space).

A~closed set $M\subset X$  is said to be  \textit{monotone path-\allowbreak con\-nec\-ted}~\cite{Ali06} if any two points
in~$M$ can be connected by a~continuous monotone curve (an~arc) $k(\cdot)\subset \nobreak M.$

That the $\mcc$-connected and monotone path-\allowbreak con\-nec\-ted sets arise natural in the study of connectedness of suns
has been demonstrated in the papers by Berens and Hetzelt and further by Brown and the author (see \cite{Brown1987}, \cite{Ali06}, \cite{Ali05}, \cite{A13FPM}):
for a~monotone path-connected  sun it proves possible to answer the long-standing question on the connectedness (acyclicity)
of the intersection of a~sun and a~ball, thereby closing the well-known Vlasov's theorem on the solarity of acyclic sets.

Given $x\in X$ and $\emptyset \ne M\subset X$,  we let $P_Mx$ denote the metric projection of~$X$ onto the set~$M$
(the set of best approximation to~$x$ from~$M$).

A subset $M\ne\emptyset$ of a normed linear space~$X$
is called a~\textit{sun}~\cite{Vla73} if, for any point,
$x\in X\setminus M$ there exists a~point $y\in P_Mx$ (a~luminosity point) such that
\begin{equation}\label{eq-1}
y\in P_M[(1-\lambda)y+\lambda x]\quad \text{for all} \quad \lambda\ge 0.
\end{equation}
The concept of a~sun was introduced by N.\,V.~Efimov and S.\,B.~Stechkin in 1958 for the study of Chebyshev sets
(note that by the term
`sun' they understood what is now called a~strict protosun---such sets~$M$ are defined by the property
that \eqref{eq-1} holds for any  $x\notin M$ and  $y\in P_Mx$).

Note~\cite {Ali12} that a~Chebyshev mono\-tone path-\allowbreak con\-nec\-ted set is always a~sun; more\-over,
if $M$~is a~mono\-tone path-\allowbreak con\-nec\-ted set and $P_Mx=\{y\}$, $x\notin M$, then $y$~is a~luminosity point.
From the definition we have that mono\-tone path-\allowbreak con\-nec\-ted\-ness is preserved if we intersect a~mono\-tone path-\allowbreak con\-nec\-ted
set with an extreme hyperplane and, moreover, with an arbitrary span (that is, a~set formed by in\-ter\-sec\-tion of hyperstrips generated
by extreme points of the dual sphere~\cite{Ali12EMJ}). A~closed ball is a~span, and hence a~mono\-tone path-\allowbreak con\-nec\-ted set is
$B$-mono\-tone path-\allowbreak con\-nec\-ted; that is, its intersection with any closed (and hence, with any open) ball is mono\-tone path-\allowbreak con\-nec\-ted.
Consequently \cite{A13FPM}, a~\textit{monotone path-\allowbreak con\-nec\-ted set is necessarily $\mcc$-connected}.
The converse assertion may fail to hold---the corresponding example is due to Franchetti and Roversi~\cite{FR}: $M=M_1\cup M_2$, where $M_i=\{x\in C[0,1]\mid x(0)=i$, $i=1,2\}$.
However, for closed sets in~$c_0$ (see~\cite {Ali05}) and for boundedly compact sets in an arbitrary separable space~$X$
these properties are equivalent  (see Theorem~1 below).

\smallskip

The present paper we shall be concerned with topological and approximative properties of $\mcc$-connected  (Menger-connected) sets and
monotone path-\allowbreak con\-nec\-ted sets and with the solarity thereof.
Theorem~\ref{t-1} asserts  that in a~broad class of Banach spaces (in particular, in separable spaces)
a~boundedly compact (boundedly weakly compact) $\mcc$-connected  set is monotone path-\allowbreak con\-nec\-ted. Furthermore,
the intersection of a~boundedly compact monotone path-connected  set with a~closed ball is cell-like (has the shape of a~point) and, in particular,
is acyclic and is a~sun.
This result partially complements the well-known Vlasov's theorem~\cite{Vla73} asserting that a~boundedly compact $P$-acyclic
subset of a~Banach space is a~sun. The case of weakly compact sets is examined in Theorem~\ref{t-2}.
In the proof of  Theorem~\ref{t-1} we will, in passing, extend the classical Rainwater--\allowbreak Simons theorem
on weak convergence of sequences to the convergence in the associated norm (Proposition~\ref{p-3.1}).

\section{Acyclic and cell-like sets}
A homology (cohomology) theory  associates with any topological space $X$ a~se\-quence of abelian groups $H_k(X)$, $k = 0, 1, 2, \dotsc$
(homology groups) and $H^k(X)$, $k = 0, 1, 2, \dotsc$ (cohomology groups), which are
homotopy invariants of a~space: if two spaces are homotopy invariant, then the corresponding (co)homology groups of
are isomorphic.
There are several ways to construct (co)homology groups of which we mention the followings: 
the construction based on nerves of covers proposed by P.\,S.~Alexandroff and further extended by E.~\v Cech;
the Vietoris construction based on the concept of true cycles and which applies to metric spaces; the construction based
on the concept of singular chains  (continuous images of simplicial chains).

Let $A$ be an nontrivial arbitrary abelian groups. A~space (all spaces are assumed to be metrizable) is called
\textit{acyclic} if its \v Cech \textit{co}homology group with coefficients from~$A$ is trivial  (it has no cycles besides boundary).
Thus, the definition of acyclicity depends on the group of coefficients.
It is worth noting that (Alexandroff--)\allowbreak \v Cech(--Dowker) homology is not a homology theory, failing to satisfy the exactness axiom,
whereas \v Cech(--Dowker) cohomology forms a~cohomology theory of to\-po\-lo\-gi\-cal spaces. For an comprehensive
account of (co)homologies on compacta, topological and uniform spaces, the reader is referred to the recent survey by Melikhov~\cite{Mel09}.

If a (co)homology has compact support  (satisfying  the compact supports axiom) and
if the coefficients of the (co)homology group lie in a~filed, then the notions of homological and cohomological acyclicity
coincide~\cite{Massey}.
However, this is not the case for an arbitrary abelian coefficient group. For example, the
2-adic solenoid (the inverse limit of the sequence
$S^1 \overset{f}\leftarrow S^1  \overset{f}\leftarrow \dotsc$, where $f=z^2$)
is acyclic  in \v Cech homology with coefficients in the field  $\mathbb Z_2$, but there is no acyclicity in \v Cech cohomology  (see, for example,~\cite{KarRe02}).

Below, unless otherwise stated,
the acyclicity will be understood in the sense of \v Cech cohomology with coefficients in an arbitrary abelian group.

A compact nonempty space is called an $R_\delta$-\textit{set} (see, for example,~\cite{Gorn}, (2.11))
if it is homeomorphic to the intersection of a decreasing sequence of compact contractible spaces (or compact absolute retracts;
see \cite{Gorn}, Theorem~2.13). $R_\delta$-sets naturally arise as spaces of solutions to the Cauchy problem and to autonomous and nonautonomous differential inclusions
\cite{Gorn00}, \cite{Dragoni},~\cite{Andres}.  Results of this kind date back to Aronszajn.

A compact space $Y$ is called  \textit{cell-like} (or having the shape of a point) if there exists an
ANR-space (absolute neighbouring retract) $Z$ and an embedding $i:Y\to Z$ such that the image $i(Y)$ is contractible
in any of its neighbourhoods  $U\subset Z$ (see \cite{Gorn}, Definition~82.4);
here, a~cell-like set need not be contractible.
From the well-known Hyman's characterization of $R_\delta$-sets it readily follows that
an $R_\delta$-set is always cell-like (\cite{HbTFT}, \S\,4.2 and \cite{Papag}, p.~50). But since any mapping of a~compact of the shape of a~point
into an ANR is homotopically trivial, we have that a~compact of the shape of a~point  (cell-like) is contractible is contractible in each of its neighbourhoods in any
ambient ANR.
As a~corollary, \textit{the classes of  $R_\delta$-sets and cell-like (having the shape of a~point) compact spaces coincide}.

Note that the cell-likeness implies the acyclicity  (with respect to any continuous  (co)homology theory; \cite{HbTFT}, p.~854),
there being acyclic sets which are not cell-like, as well as cell-like sets which are not path-connected (topologist's sine curve).

\smallskip\goodbreak

Following Vlasov~\cite{Vla73}, if $\textrm{Q}$ denotes some property
(for example, `connected'), then we say that  $M$~is
\begin {list}{}{\leftmargin=15pt \itemsep=0pt\topsep=2pt\parsep=0pt}
\item[] $P$-Q if $P_M(x)$ is nonempty and is~$\textrm{Q}$ for all $x\in X$;
\item[] $B$-Q if  $M\cap B(x,r)$ is~Q for all $x\in X$, $r>0$;
\item[] $\mathring B$-Q if  $M\cap \mathring B(x,r)$ is~Q for all $x\in X$, $r>0$.
\end{list}
Thus a closed subset of a~finite-dimensional space is $P$-nonempty, or is an existence set (proximinal).
Correspondingly, a~set $M$ is $P$-\textit{acyc\-lic} if $P_Mx$ is nonempty and is acyclic for any~$x$;
a~set $M$ is $B$-\textit{acyc\-lic} if the intersection of~$M$ with any closed ball is acyclic.

\begin{remark} \rm $B$-connected sets are also called $V$-connected set (here the letter `V' comes from the works of L.\,P.~Vlasov,
who denoted balls by  $B(x,r)$). Our term $B$-connectedness agrees with the more conventional notation in which balls are denoted by $B(x,r)$,
as well as with the term `bounded connectedness', which was  introduced by D.~Wulbert in the 1960th.
\end{remark}

Brown \cite{Brown1986}, Corollary 1.6.2, proved that if a boundedly compact subset $M$ of a~Banach space
is $P$-acyc\-lic   (relative to \v Cech cohomology with coefficients in an arbitrary abelian group),
then $M$ is $B$-acyc\-lic\footnote{Actually, Brown \cite{Brown1986} has proved slightly more:
if $M$ is a~$P$-acyclic approximatively compact subset of a~Banach space and
if the intersection of~$M$ with some ball $B$ is compact, then  $M\cap B$ is acyclic.}.
Consequently, the acyclicity of an arbitrary compact
$\mcc$-connected  set implies the $P$- and $B$-acyc\-lic\-ity
of an arbitrary boundedly compact $\mcc$-connected  set~$M$. This follows from that fact that the intersection
of such an~$M$ with an arbitrary ball  $B(x, r)$ is compact and $\mcc$-connected.

\section{Banach--Mazur hull. Monotone path-\allowbreak con\-nec\-ted\-ness. Spaces of classes $(\stone)$ and $(\sttwo)$.
Associated norm. Rainwater--\allowbreak Simons's theorem}

We examine the relationship between the Banach--Mazur hull $\mcc(\cdot,\cdot)$ and function intervals
$\brL \cdot,\cdot\brR$ (to be defined in~\eqref{eq-3} below). We introduce the class $(\stone)$ of normed linear spaces (including all separable normed linear spaces)
in which $\mcc(\cdot,\cdot) =\brL \cdot,\cdot\brR$.
An important property of the Banach--Mazur hull $\mcc(\cdot,\cdot)$ is that
$z\in \mcc(x,y)$ if and only if $z$ lies metrically between $x$ and~$y$ with respect to the so-called (Brown)-associated norm~$|\cdot|$.
Such a~norm is introduced on the spaces from  the class $(\sttwo)$, which also includes all separable normed linear spaces.
This observation enables us to use the machinery of metrical convexity (Lemma~\ref{l-5.A}).
In passing, we generalize the classical Rainwater--\allowbreak Simons  theorem on weak convergence of sequences to
the convergence with respect to the associated norm (Proposition~\ref{p-3.1}).

To begin with, we note that in  $X=C(Q)$ the structure of the hull $\mcc(x,y)$, $x,y\in X$, is quite apparent:
\begin{equation*}
  \mcc(x,y)  =\brL x,y\brR,
\end{equation*}
where
\begin{equation}\label{eq-2}
  \brL x, y \brR: =\bigl\{z\in C(Q) \mid z(t)\in [x (t), y(t)]
  \quad \forall \, t\in Q\bigr\}.
\end{equation}
A similar representation also holds in  $C_0(Q)$-spaces
($Q$ is a locally compact Hausdorff space)---this
follows from a characterization of extreme elements of the unit sphere of the dual
space to $C_0(Q)$, which is due to B.~Brosowski, F.~Deutsch and P.\,D.~Morris~\cite{BD74}.

In the general case, the interval  $\brL x,y\brR$ is defined, in analogy with~\eqref{eq-2}, as follows:
\begin{equation}\label{eq-3}
\brL x,y\brR := \{z\in X\mid \min \{\varphi(x),\varphi(y)\} \le
\varphi(z)\le  \max \{\varphi(x),\varphi(y)\} \quad
  \forall \varphi \in \operatorname{ext} S^*\}.
\end{equation}

Following Brown \cite{Brown1987} and Franchetti and Roversi \cite{FR} we consider the
class (\stone) of spaces~$X$ such that
$$
   \mcc(x,y)=\brL x,y\brR\quad \text{for all}\quad  x,y\in X.
\leqno {(\stone)}
$$

The abbreviation $(\stone)$ comes from the phrase `The hull $\mcc(x,y)$ equals the interval $\brL x,y\brR$ for all $x,y$'. The author does
not know examples of Banach spaces which are not $(\stone)$-spaces.

Note that the inclusion $ \mcc(x,y)\supset \brL x,y\brR $
holds in any normed linear space~$X$ (see, for example, \cite{A13FPM}).
It is also known \cite{A13FPM}, \cite{FR} that the equality  $\mcc(x,y)=\brL x,y\brR$ holds in a~broad class of Banach spaces,
and in particular, in spaces for which smooth points are dense on the unit sphere
(such a~class includes the weakly Asplund spaces, and hence, the separable spaces).

It is readily verified that the class $(\stone)$ contains the spaces $C(Q)$ on  a~Haus\-dorff compact set~$Q$,
and in particular, the space~$\ell^\infty$ (\textit{qua} the space of continuous functions
on the Stone--\v Cech compactification of the natural numbers).
We also note that if a~space~$X$ is such that  $\operatorname{ext} S^*$ lies in the closure of the set of  $w^*$-semi-denting
points of the dual ball~$B^*$ (\textit{Moreno's condition}), then  $\brL x,y\brR = \mcc(x,y)$ for all $x,y\in X$;
this condition, for example, is satisfied for finite-dimensional spaces and spaces with the Mazur intersection property.
Recall that a~point ~$f \in S^*$ is called a~$w^*$-\textit{semi-denting} point of the dual ball~$B^*$ (see, for example, \cite{Giles00}) if, for any
 $\varepsilon >0$, there exists a~$w^*$-slice~$S\ell$ of~$B^*$ such that $\operatorname{diam}(\{f\}\cup S\ell)< \varepsilon$.
Here,  $S\ell(B^*,x,\delta):= \{g\in S^*\mid g(x)> 1-\delta\}$, $0<\delta<1$, $x\in X$.

We shall also require the following class of spaces introduced by Franchetti and Roversi~\cite{FR}:
$$
   \operatorname{ext} S^* \ \ \text{is} \ \ w^*\text{separable}.\leqno (\sttwo)
$$
In the definition of class $(\sttwo)$ we will always assume that
$$
  F=(f_i)_{i\in I}\subset \operatorname{ext} S^* \ \ \text{is}\ \ w^*\text{-dense in}\ \operatorname{ext} S^*,\quad
  \operatorname{card} I\le \aleph_0,\quad F= -F.
$$

The abbreviation $(\sttwo)$ is taken from the German  `Die Extrem\-punkt\-menge der konjugierten Einheitskugel ist $w^*$-separabel'.
\smallskip

From Krein--Milman theorem it readily follows that any space from the class  $(\sttwo)$ has  $w^*$-separable dual ball  $B^*$.
The latter condition is equivalent \cite{Dancer} to the fact that $X$~is isometrically isomorphic to a~subspace of~$\ell^\infty$.
We do not know a~characterization of the class $(\sttwo)$.
We also point out  that there are examples of  $C(K)$-spaces, $K$~is a~nonseparable Hausdorff compact set or
spaces of the form $X=\ell_1\oplus \ell_2(\Gamma)$, $|\Gamma|=\mathfrak c$, for which  $X^*$ is
$w^*$-separable, but the dual unit ball~$B^*$ is not.

Further, it is well known that if  $X$~is a~separable normed linear space, then the $w^*$-topology of the
dual unit ball $B^*$ is metrizable. As a~result, the ball $B^*$ is $w^*$-separable (\cite{Fabian}, Corollary 3.104).
Hence, \textit{any separable space lies the class}$(\sttwo)$. The
class  $(\sttwo)$ also contains the \textit{nonseparable} space~$\ell^\infty$ (\textit{qua} the space of continuous functions on
the Stone--\v Cech compactification  $\beta\mathbb N$ of the natural numbers--- such compact set is separable but nonmetrizable).
Also note that the $C(Q)$ on a~nonseparable~$Q$ and $c_0(\Gamma)$ on an uncountable~$\Gamma$ fail to lie in~($\sttwo$).

It would be interesting to characterize the class  $(\sttwo)$.

In summary, with regard to spaces of the class  $(\stone)$ and $(\sttwo)$, we point out that
$$
 \parbox{100mm}
 {the class $(\stone)\cap (\sttwo)$ contains all separable Banach spaces
(and in particular,  $C(Q)$ on a~metrizable compact set~$Q$) and the nonseparable space~$\ell^\infty$).}
$$

Let  $X\in (\sttwo)$, let $F=(f_i)_{i\in I}$ be the family of functionals from the definition of the class~$(\sttwo)$, let
$(\alpha_i)\subset \mathbb{R}$, $\alpha_i>0$, $i\in I$, and let
$\sum \alpha_i <\infty$. Given $x\in X$, we set
\begin{equation}\label{eq-4}
 |x| = \sum _{i\in I}\alpha_i |f_i(x)|.
\end{equation}
It is easily seen that $|{\,\cdot\,}|$ defines a norm on~$X$, which we shall call, following
Brown~\cite{Brown1987}, the \textit{associated} norm on~$X$. Clearly,  $|x|\le \|x\|\sum \alpha_i$.

The importance of the associated norm can be seen from the following result, which is a direct and straightforward generalization
of Corollary 3.2 of~\cite{Brown1987}, which was put forward by A.\,L.\,Brown for finite-dimensional~$X$.

\begin{lemma} \label{l-3.1}
Let $X\in (\stone) \cap (\sttwo)$ be a~Banach space {\rm (}in particular, $X$~is a~separable
Banach space{\rm )} and let $x,y\in X$. Then the following conditions are equivalent{\rm :}
\begin {list}{}{\itemsep=0pt\topsep=2pt\parsep=0pt}
\item[\rm  a)] $z\in \mcc(x,y);$
\item[\rm  b)] $|f_i(x)-f_i(y)|=|f_i(x)-f_i(z)|+|f_i(z)-f_i(y)|$ for all
    $i\in I$, where $F=(f_i)_{i\in I}$ is the family from the definition of $(\sttwo);$
\item[\rm  c)] $|x-y|=|x-z|+|z-y|$ {\rm (}that is,  $z$ is $|{\,\cdot\,}|$-between $x$ and~$y)$.
\end{list}
\end{lemma}

The following result may be regarded as an extension of the classical
Rain\-water--\allowbreak Simons theorem (see, for example, \S\,3.11.8.5 in~\cite{Fabian})
to the convergence in the associated norm~${|{\,\cdot\,}|}$ on spaces from the class $(\sttwo)$ (in particular, on separable normed linear spaces).
Rainwater--\allowbreak Simons's theorem  states that
states that a~bounded sequence $(x_n)$ in a~Banach space~$X$
weakly converges to an  $x \in X$ if and only if the sequence $(f (x_n))$ converges to~$f (x)$ for any functional~$f$
from an arbitrary fixed James boundary for~$X$ (for example, for all $f\in \operatorname{ext} S^*$).
Thus, even though the weak convergence is nonmetrizable in general, there is a~norm on~$X\in (\sttwo)$ with respect to which
the convergence of  \textit{sequences} is equivalent to the weak convergence.

Here, we recall (see, for example, \S\,3.11.8 in~\cite{Fabian}) that
a~subset~$A$ of the dual unit sphere $S^*$ of~$X^*$
is called a~(James) \textit{boundary} for the space~$X$ if, for any $x\in X$, there exists an  $f\in A$ such that $f(x)=\|x\|$.
It is an easy consequence of Krein--\allowbreak Milman's theorem that the set $\operatorname{ext} S^*$ of extreme points of the dual unit ball is a~boundary for~$X$.

\begin{propos}\label{p-3.1} Let $X\in (\sttwo)$ be a Banach space, $F:= (f_i)_{i\in I}\subset \operatorname{ext} S^*$ be the family of functionals
from the definition of $(\sttwo)$. Also let $(x_n)$ be a~bounded subsequence in~$X$.
Consider the following conditions\/{\rm :}
\begin {list}{}{\itemsep=0pt\topsep=2pt\parsep=0pt}
\item [{\rm a)}] $x_n\overset{|{\,\cdot\,}|}{\longrightarrow}x$\/{\rm ;}
\item [{\rm b)}]  $f_i(x_n)\to f_i(x)$ for any~$i\in I$\/{\rm ;}
\item [{\rm c)}] $x_n\overset{w}{\longrightarrow}x$.
\end{list}

Then conditions {\rm a)} and {\rm b)} are equivalent, either of which follows from~{\rm c)}.
If  $X^*$ is separable, then all three conditions are equivalent.
\end{propos}

\begin{remark} \rm That $X^*$ is separable in  b)$\Rightarrow$c) is essential. Indeed, let  $X=\ell^1$.
Consider finite sequences from~$\ell^\infty$ consisting of zeros and $\pm 1$. This set is countable and is $w^*$-dense in $\operatorname{ext} S^*$.
However, the convergence of elements from~$\ell^1$ on these sequences does not imply their weak convergence. This fact was noted by P.\,A.~Bor\-odin
in disputing the results of the paper.
\end{remark}

\begin{remark} \rm Proposition \ref{p-3.1} implies that $X,|\,\cdot\,|)$ is always a~Schur space
with respect to the associated norm $|\,\cdot\,|$  (recall that a~space is a~\textit{Schur space} if weakly convergent sequences in X are norm convergent;
the space~$\ell^1$ is a~classical example of a~Schur space).
Indeed, from the definition of the associated norm it follows that $|x|>\alpha_i |f_i(x)|$ for any~$i$. Hence, any  $f_i\in F$ lies in $X^*_{|\,\cdot\,|}$. Now
it follows from assertion~b) of Proposition \ref{p-3.1} that if ($x_n)$\enskip
$w_{|\,\cdot\,|}$-converges, then  ($x_n)$\enskip $|\,\cdot\,|$-converges. Thus, for a~bounded sequence ($x_n)$
$$
  x_n\overset {w_{|\,\cdot\,|}}{\longrightarrow} x \iff x_n\overset {|\,\cdot\,|}{\longrightarrow} x
$$
As a result, in $(X,|\,\cdot\,|)$\enskip $w_{|\,\cdot\,|}$-compactness coincides with the strong $|\,\cdot\,|$-compact\-ness. Thus, a~reflexive  $(X,|\,\cdot\,|)$ is
finite-dimensional. This note note arose during conversations between the author and O.~Nygaard, to whom the author wishes to express his thanks.
\end{remark}

\noindent \textbf{Proof  of Proposition~\ref{p-3.1}.}
\rm The implication a)$\Rightarrow$b) is quite clear: if $x_n\overset{|{\,\cdot\,}|}{\longrightarrow}x$
(in the norm $|{\,\cdot\,}|$), then the sum
$\sum \alpha_i |f_i(x_n)-f_i(x)|$ is small for all sufficiently large~$n$; as a~corollary,
for each fixed~$j$ the difference $|f_j(x_n)-f_j(x)|$ is also small for such~$n$.

Let us prove b)$\Rightarrow$a). For any $n$, in the sum
\begin{equation}\label{eq-5}
\sum_{i\in I} \alpha_i |f_i(x_n)-f_i(x)|
\end{equation}
we set $a_i(n) := \alpha_i |f_i(x_n)-f_i(x)|$ and split the sum \eqref{eq-5} into two: for  $i\le N$ and for  $i>N$ ($N$~will be chosen later).
By the hypothesis, the sequence $(x_n)$ is uniformly bounded, and hence, in the second sum we have $f_i(x_n)-f_i(x) |\le C$ (where $C$~is independent of~$i, n$).
So, the second sum is bounded from above by the sum $\sum _{i>N} C\alpha _i<\infty$. From $\varepsilon>0$ we choose an~$N$ for which the second sum
is smaller than~$\varepsilon$. The first sum is finite, and there we choose large~$n$.

That c)$\Rightarrow$b) is clear. Assume that $X^*$ is separable and prove that b)$\Rightarrow$c).
Note that the following assertions are equivalent for a~Banach space~$X${\rm :}
 \begin {list}{}{\itemsep=0pt\topsep=2pt\parsep=0pt}
\item [-] $X$ has a~separable boundary{\rm ;}
\item [-] the boundary $\operatorname{ext}B^*$ is separable{\rm ;}
\item [-] space $X^*$ is separable.
\end{list}

Further, we shall need the concept of a  (I)-generating set introduced by Fonf and Lindenstrauss.
By definition, a~set~$C\subset B^*$ \enskip (I)-\textit{generates} the dual ball $B^*$ if
\begin{equation}\label{eq-6}
  B^* = \overline{{\rm conv}\vphantom{f}}\Bigl( \bigcup\nolimits_i \overline{{\rm conv}\vphantom{f}}\,^{w^*}C_i\Bigr)
\end{equation}
for any representation $C=\bigcup C_i$ as a~countable union of sets~$C_i$. In this definition, `I' comes from the Latin \textit{inter\-medi\-us}
and is explained by the fact that
\begin{equation*}
    B^* = \overline{{\rm conv}\vphantom{f}} C \implies C\enskip({\rm I})\text{-generates}\ B^* \implies  B^*=\overline{{\rm conv}\vphantom{f}}\,^{w^*} C.
\end{equation*}

Let $C_i:=\{f_1,\dots, f_i\}$, $i\in I$ ($F:= (f_i)_{i\in I}$ is the family of functionals from the definition of the class  $(\sttwo)$). Clearly, $F=\bigcup C_i$.
By (Kadets--Fonf--) Gode\-froy--\allowbreak Rod\'e theorem,
\begin{equation}\label{eq-7}
B^*=\overline{{\rm conv}\vphantom{f}}\,^{\|\cdot\|}\operatorname{ext} B^*.
\end{equation}
The space $X^*$  is separable,  and hence so is  $\operatorname{ext}B^*$. By the definition of the class $(\sttwo)$,
$\operatorname{ext} S^* \subset     \overline{F}{}^{w^*}$,
and  hence, $\operatorname{ext} S^* \subset \bigl( \bigcup\nolimits_i \overline{{\rm conv}\vphantom{f}}\,^{w^*}C_i\bigr)$. Finally,
$F$\enskip  (I)-gene\-rates the ball $B^*$ by \eqref{eq-7} and~\eqref{eq-6}.


Now it remains to invoke one result, obtained independently by Nygaard \cite{Nygaard} and Kalenda \cite{Kal}, which asserts that if
$C$\enskip  (I)-gene\-rates the ball $B^*$, then $C$ is a~Rainwater set; that is, a~set
with the property: if a~bounded sequence $(x_n)\subset X$ converges pointwise on~$B$, then $(x_n)$ converges weakly.
The proof of Proposition~\ref{p-3.1} is complete.
\hfill$\Box$

\section {Statements of the main results}

The mains results of the paper are given are follows.

\begin{theorem} \label{t-1}
Let $X \in (\stone) \cap (\sttwo)$ be a Banach space {\rm (}in particular, $X$~is a~separable
Banach space{\rm )} and let
set  $M\subset X$ be closed and $\mcc$-connected. Assume that at least one of the following conditions
is satisfied{\rm :}
\begin {list}{}{\itemsep=0pt\topsep=2pt\parsep=0pt}
\item[\rm  a)] $M$ is boundedly compact  {\rm (}in the norm $\|\cdot\|$ of $X);$
\item[\rm  b)] $M$ is $|{\,\cdot\,}|$-closed and
$\mcc(x,y)$ $|{\,\cdot\,}|$-compact for all $x,y\in X${\rm ;}
\item[\rm  c)] $\mcc(x,y)$ is
$\|\cdot\|$-compact for all  $x,y\in X${\rm .}
\end {list}
Then $M$ is monotone path-\allowbreak con\-nec\-ted, $P$- and $B$-cell-like,
$P$- and $B$-acyc\-lic {\rm (}relative to any continuous
{\rm (}co{\rm )}homology theory{\rm )} and is a~sun.

If $X$ is finite-dimensional, then in addition $M$ is  $P$- and $B$-contractible.
\end{theorem}

Recall \cite{TsarEC} that a set $M\subset X$ has continuous multiplicative  (additive) $\varepsilon $-selection
of the metric projection for all $\varepsilon > 0$ if and and only if $M$ is  $\mathring  B$-contractible ($B$-contractible, if $M$ is approximatively compact).
Now from Theorem~\ref{t-1} it follows that the metric projection onto an $\mcc$-connected (monotone path-\allowbreak con\-nec\-ted)
closed subset of a~finite-dimensional space has continuous multiplicative (additive) $\varepsilon $-selec\-tion for all
$\varepsilon > 0$. There is some hope that the similar result also holds in the in finite-dimensional setting.

For spaces with linear ball embedding  \cite{Ali12EMJ} (in particular,
in the spaces $\ell^1(n)$, $C(Q)$, $Q$~is a~metrizable compact set and $C_0(Q)$) Theorem~\ref{t-1} partially extends
the results of Balashov and Ivanov (\cite{BalIva09}, Theorem 2.9 and Lemma  4.18) on the path-connectedness of $R$-weakly convex (in the sense of Vial)
sets:  in a~space with linear ball embedding the intersection  of an~$R$-weakly convex set with a~closed or
open ball is $\mcc$-connected (\cite{Ali12EMJ}, Theorem 4.1).

Theorem~\ref{t-1} has the following corollary.

\begin{corollary}
Let $X \in (\stone) \cap (\sttwo)$ be a Banach space {\rm (}in particular, $X$~is a~separable
Banach space{\rm )} and let
set  $M\subset X$ be closed. Assume that for some $x\notin M$ \enskip  $P_Mx$ is compact and $\mcc$-connected.
Then $M$ is monotone path-\allowbreak con\-nec\-ted, $P$- and $B$-cell-like,
$P$- and $B$-acyc\-lic and is a~sun.
\end{corollary}

In this connection, we note that
$$
\text{if}\ M\subset X\in (\stone) \ \text{is $\mcc$-connected, then so is}\ P_Mx.
$$

Indeed, let  $u,v\in P_Mx$. By the assumption  $M$ is  $\mcc$-connected, and so there exist $ z\in \mcc(u,v)\cap M$, $z\ne u,v$.
Assume that  $z\notin S(x,\|x-u\|)$. Since  $\operatorname{ext} S^*$ is a~boundary for~$X$,
any point outside a~closed ball (the ball $B(x,\|x-u\|)$ in our setting) can be strictly separated from it by an extreme functional. Hence, $z\notin\brL u,v\brR$.
This is a~contradiction, because  $z\in \mcc(u,v)$, while in $X\in (\stone)$ one always has $\mcc(\cdot,\cdot) = \brL\cdot,\cdot\brR$.

\smallskip \goodbreak

For weakly compact sets our result is substantially weaker.

\begin{theorem} \label{t-2} Let $X$~be a separable Banach space and let  $\emptyset\ne M\subset X$ be boundedly weakly compact.
Assume that $M$ is $\mcc$-connected. Then $M$ is monotone path-\allowbreak con\-nec\-ted.
\end{theorem}

\begin{problem} \rm The following result due to Vlasov's is well-known~\cite{Vla73}: \textit{a~$P$-acyc\-lic
boundedly compact subset of a~Banach space is sun}\footnote
{The proof of Vlasov's theorem depends on the classical
Eilenberg--\allowbreak Montgomery's fixed point theorem (see, for example, \cite {Gorn}, Corollary (32.12)) in which
the acyclicity is understood in the sense of (Alexandroff--)\allowbreak \v Cech homology groups with coefficients in a~field
(in the paper by Eilenberg--\allowbreak Montgomery the acyclicity is understood in the sense of Vietoris cycles and
homology groups over a~field). However, \v Cech homology groups are known to be isomorphic to Vietoris homology groups on the
category of compact metrizable spaces, and the notions of the homological and cohomological acyclicity coincide if the
coefficients of the homology (cohomology) groups lie in a~field and if the support is compact.}.
It will be interesting, especially in view of Theorem~\ref{t-2}, to extend
Vlasov's theorem to the case of boundedly weakly compact sets, at least in the separable setting.  The difficulty here is as follows:
despite the fact that a~weakly compact subset~$M$ of a~separable Banach space is metrizable and the topology
generated on~$M$ by the metric agrees with the weak topology on~$M$, Eilenberg--\allowbreak Montgomery's theorem
does not apply, because the metric projection to~$M$ is only norm-to-weak upper semicontinuous  (the support is the same, but the topologies on it are different).
We also note in passing that the problem of the solarity of weakly compact Chebyshev sets.
has not yet been solved in the general setting.
\end{problem}

\section {The proofs}

The idea of Theorem~\ref{t-1} goes back to~\cite{FR}.
The monotone path-\allowbreak con\-nec\-ted\-ness is established in~\cite{Ali06}, and the remaining assertions follow from Lemma~\ref{l-5.2},
which will be proved below. The finite-dimensional case is a~consequence of assertion~a) and Brown's theorem~\cite{Brown1988}, in accordance to which
the intersection of an $\mcc$-connected closed subset of a~finite-dimensional~$X$  with a~closed ball is $n$-connected for all $n\in \mathbb Z_+$ (that is, each map into the set
of a~$k$-sphere, $k\le n$, extends continuously to the $k+1$ ball).
Hence, according to the well-known characterization of ARs (see, for example, \cite{Hu}, Theorem 11.1), $M$ is contractible and locally contractible.
Note~\cite{Ali05} that condition~b) of Theorem~\ref{t-1} holds \textit{a~fortiori} in the space $X=c_0$.

\goodbreak

We shall need the following auxiliary result.

Recall that a set $M$ is called \textit{metrically convex} or  $d$-\textit{convex} (Menger-convex)  \cite{Men}
with respect to metric~$d$ if, for any distinct points  $x,y\in M$, the set $M\setminus\{x,y\}$ contains the point~$z$ lying
 $d$-between $x$ and~$y$; that is,  $d(x,y) = d(x,z)+d(z,y)$. The following result was established by Menger \cite{Men}, see also \cite{Kirk}, p.~24.

\begin{lemma} \label{l-5.A} Let $(Y,d)$ be a~complete Menger-\allowbreak  convex metric space.
Then for each  $x$ and $y$ in~$Y$ there exists an isometry $f:[0, d(x,y)]\to Y$ such that
$f(0)=x$ and $f(d(x,y))=y$. Consequently, $Y$ is path-connected.
\end{lemma}

\noindent \textbf{Proof of Theorem~\ref{t-2}.} \rm It suffices to deal with  a~bounded set~$M$.
It is well known  (see, for example, \cite [Pro\-po\-si\-tion 3.107]{Fabian}) that
if $M$ is a~weakly compact subset of a~Banach space with $w^*$-separable dual~$X^*$,
then $M$ with the relative weak topology is metrizable. It is also known that
a~weakly compact subset of a~Banach space is weakly complete  (sequentially weakly complete).
In view of Proposition~\ref{p-3.1}, in spaces of class $(\sttwo)$, and hence, in any separable space~\cite{A13FPM},
$w$-completeness implies the $|\cdot|$-completeness with respect to the (Brown)-associated norm~$|\cdot|$.
Now in order to apply Lemma~\ref{l-5.A} it remains to observe that in a~separable space the $\mcc$-connectedness of a~set
is equivalent to its $|\cdot|$-convexity (see Lemma~\ref{l-3.1}).

Thus, according to Lemma \ref {l-5.A} any two points from~$M$ are connected by an arc which is  the range of an isometry. We claim that
this arc is monotone.

Let $x,y\in M$. We need to show that $x$ and $y$ are connected by monotone curve in~$M$.
Without loss of generality we assume that  $x=0$, $f(y)\ge 0$. Let $k(\cdot): [0,|x-y|]\to X$, $k(0)=x$, $k(|x-y|)=y$, be an isometry existing
by Lemma~\ref{l-5.A}. Further, by Lemma~\ref{l-5.A} for each  $\delta\in
(0,1)$ there exists a~point  $u\in k(\cdot)$ such that $|u|=\delta$,
$|y-u|=\delta$. Since $k(\cdot)$ is an isometry, we have
\begin{equation*}
  |y|=\alpha_1 |f(u)| + \alpha_1 |f(y)-f(u)| +
  \sum _{i\ge 2} \alpha_i \bigl(|f_i(u)|+|f_i(y)-f_i(u)|\bigr).
\end{equation*}
Elementary calculations show that
$$
  |y|=
  \alpha_1 |f(u)| + \alpha_1 |f(y)-f(u)| \ge
  \alpha_1 |f(u)| + \alpha_1 |f(y)-f(u)| +
  \sum _{i\ge 2}\alpha _i|f_i(y)|.
$$
Therefore assuming that $f(u)<0$ or $f(u)> f(y)$ we obtain
$|y|>|y|$, a~contradiction. Thus, for any $f\in \operatorname{ext} S^*$, the function  $f(k(t))$ is monotone in~$t$, as required.

Let  $X \in (\stone)$ and $(\sttwo)$ and
let $F=(f_i)\subset \operatorname{ext} S^*$ be a~$w^*$-dense in $\operatorname{ext} S^*$ family from the definition of the class $(\sttwo)$. Given any $n\in \mathbb{N}$,
we consider a~bounded linear operator $s_n:X\to \ell^\infty(n)$ defined as
$$
  s_n(x)=\bigl(f_1(x),\dots,f_n(x)\bigr).
$$
Note that $\|s_n(x)\|\le \|x\|$ and $\|s_n(x)\|\to \|x\|$ as  $n\to \infty$, since  $\exp S^*$ is a~James boundary for~$X$;
that is,  $\|x\|=\max \{f(x)\mid f\in \operatorname{ext} S^*\}$.

The first assertion in Lemma~\ref{l-5.1} is proved by a slight modification of the argument given by Franchetti and Roversi~\cite{FR}.

\begin{lemma} \label{l-5.1}
Let $X\in  (\stone) \cap (\sttwo)$ be a~Banach space {\rm (}in particular, $X$~is a~separable
Banach space{\rm )} and let $\emptyset\ne M\subset X$ be boundedly compact. Then the following results hold:
\begin {list}{}{} 
\item[\rm 1)] if $M$ is $\mcc$-connected in~$X$, then $s_n(M)$ is
monotone path-\allowbreak con\-nec\-ted in $\ell^\infty(n)$ for all $n\in \mathbb{N}${\rm ;}
\item[\rm 2)] if $s_n(M)$ is $\mcc$-connected in $\ell^\infty(n)$ for all $n\in \mathbb{N}$, then
$M$ is monotone path-\allowbreak con\-nec\-ted in~$X$.
\end{list}\end{lemma}

Note that assertion~2) of Lemma~\ref{l-5.1} is not used in the proof of Theorem~\ref{t-1} and Lemma~\ref {l-5.2}.

\begin{lemma} \label {l-5.2} Let  $X\in (\stone)\cap (\sttwo)$ be a~Banach space
{\rm (}in particular, $X$ separable or  $X=\ell^\infty)$
and let a~set $\emptyset\ne M\subset X$ be boundedly compact and $\mcc$-connected. Then
$M$\enskip $P$- and $B$-cell-like; that is, each of the sets
$$
  P_Mx, \ \  M\cap B(x,r),  \ \ x\in X,  \ \ r>0,
$$
is cell-like. In particular, the set $M$\enskip $P$-
and $B$-acyclic {\rm (}with respect to any continuous  {\rm (}co{\rm )}homology theory{\rm )}.
\end{lemma}

Before proceeding with the proof of Lemma~\ref{l-5.2} we recall that
an  \textit{inverse system {\rm (}inverse spectrum  {\rm if  the system is countable)} of topological spaces} (see, for example, \cite{EGT}, p.~56)
is the family  $\mathcal S=\{X_\alpha, \pi^\beta_\alpha,
\Sigma\}$, where the set $\Sigma$ is partially ordered by~$\prec$,
$X_\alpha$ is a~topological Hausdorff space, and
$\pi^\beta_\alpha: X_\beta \to X_\alpha $ is a~continuous mapping for
any $\alpha\prec \beta$ such that
$\pi^\alpha_\alpha=\operatorname {id} _{X_\alpha}$ and
$\pi^\beta_\alpha\pi^\gamma_\beta= \pi^\gamma_\alpha$ for all $\alpha \prec
\beta \prec\gamma$.
A subspace of $\prod_{\alpha\in \Sigma} X_\alpha$ is called the inverse limit of~$\mathcal S$:
$$
  \mathop {\vphantom{}
  \underleftarrow{\lim}} \mathcal S=
  \biggl\{ (x_\alpha)\in \prod_{\alpha\in \Sigma} X_\alpha \ \bigg|\
  \pi^\beta_\alpha(x_\beta) =x_\alpha\ \ \text{for all}\ \ \alpha\prec\beta
  \biggr\}
$$

If  $\pi_\alpha:\underleftarrow{\lim} \mathcal S\to X_\alpha$ is the restriction of the projection
$p_\alpha:\prod_{\alpha\in\Sigma} X_\alpha\to X_\alpha$ onto the $\alpha$th axis  (the canonical projection), then
$\pi_\alpha=\pi_\alpha^\beta \pi_\beta$ for all $\alpha\prec\beta$,

We shall need the following result, in which assertion 1) is contained, for example, in~\cite{Dragoni}, and assertion 2) is obtained in~\cite{Gabor} (see also~\cite{Andres}).
Also note that assertion~a) is a corollary to the classical Tikhonov's compactness theorem, and
assertion~b) is well known and belongs to A.\,G.~Kurosh and N.\,E.~Steenrod.

\begin{lemma} \label{l-5.B}
{\rm 1)} Let $\mathcal S=\{X_\alpha, \pi^\beta_\alpha, \Sigma\}$ be an inverse system. Then the limit $\mathop {\vphantom{}
\underleftarrow{\lim}} \mathcal S$ is a~closed subset of the product  $\prod_{\alpha\in \Sigma} X_\alpha$. Moreover, if, for every  $\alpha\in \Sigma$,
\begin {list}{}{\itemsep=0pt\topsep=2pt\parsep=0pt}
\item[{\rm a)}] $X_\alpha$ is compact, then so is the limit  $\mathop {\vphantom{}
\underleftarrow{\lim}} \mathcal S$;
\item[{\rm b)}] $X_\alpha$ is compact and nonempty, then so is  $\mathop {\vphantom{}
\underleftarrow{\lim}} \mathcal S$;
\item[{\rm c)}] $X_\alpha$ is a continuum,  then so is $\mathop {\vphantom{}
\underleftarrow{\lim}} \mathcal \mathcal S$;
\item[{\rm d)}] $X_\alpha$ is compact and acyclic\footnote{With respect to any continuous theory of cohomology.},
then so is  $\mathop {\vphantom{} \underleftarrow{\lim}} \mathcal S$;
\item[{\rm e)}] $X_\alpha$ is metrizable and $\Sigma$ is countable, then the limit $\mathop {\vphantom{} \underleftarrow{\lim}} \mathcal S$ is metrizable.
\end{list}
{\rm 2)} If $\mathcal S=\{X_n, \pi^p_n, \mathbb{N}\}$ is an inverse spectrum and each $X_n$
is compact and cell-like, then the inverse limit
$\mathop {\vphantom{}
\underleftarrow{\lim}} \mathcal S$ is compact and cell-like.
\end{lemma}

\noindent \textbf{Proof  of Lemma~\ref{l-5.2}.}
It is sufficient to consider the case when $M$ is compact.
 Let $\Omega$~be the family of finite subsets from $F\subset \operatorname{ext} S^*$.
The set of all finite sequences of natural numbers is countable, and so is~$\Omega$.
The set $\Omega$ is directed by inclusion; that is, by the relation~$\prec$ defined by letting  $A\prec B$ if and only if $A\subset B$
Given  $A\in  \Omega$, we define a~map  $s_A:X\to \ell^\infty(A)$ by setting
$$
 s_A(x) =\{f(x)\mid f\in A\}.
$$
Then $s_A$ is a~bounded linear operator.
Given $A, B \in \Omega$, $A\supset B$, there is the natural restriction mapping  $\ell^\infty(A)\to \ell^\infty(B)$.
The family $\{s_A\mid A\in \Omega\}$ forms an inverse system relative to the restriction mappings. By assertion~1) of Lemma~\ref{l-5.1}, each  $s_A(M)$ is
compact $\mcc$-connected subset of the space~$\ell^\infty(A)$, and further, by
Brown's theorem (\cite{Brown1986}, Theorem~1), is infinitely connected (and hence, is cell-like and acyc\-lic).

Since $\{s_A\mid A\in \Omega\}$~is an inverse spectrum, the triangular diagram
$ s_A =s_{BA}\circ s_B$ is commutative for  $A,B\in \Omega$, $A\subset B$, where $s_{BA}$~is the restriction mapping (see  \cite[p.~428]{Dug}).
Hence  $g:M\to M_\infty:= \mathop {\vphantom{} \underleftarrow{\lim}} s_A(M)$ is continuous (by the above) and is injective  (see \cite[(6.9)]{Dug}).

We claim that $g$ is surjective. Assume that
$$
  (x_A)_{A\in \Omega} \in
  \mathop {\vphantom{} \underleftarrow{\lim}} s_A(M)
  \subset \prod _{A\in \Omega} s_A(M).
$$
For each $A\in M$ we choose $y_A\in M$ such that  $s_A(y_A)=x_A$.
Then $(y_A)_{A\in\Omega}$ is a~net in~$M$, and hence has a~cluster point $y\in M$, since $M$ is compact. It remains
to prove that $s_A(y)=x_A$ for any $A\in \Omega$. Let $A,B\in \Omega$, $B\supset A$. Then $s_A(y_B)=s_B(y_B)|_A=x_A$. Passing to the limit over a~subnet we see that $s_A(y)=x_A$.

Thus, $g:M\to M_\infty$ is continuous and is a bijection.
Hence, $M$ is homeomorphic to~$ M_\infty$.

Further, cell-likeness (as well as acyclicity) is a topological property. We have $M\simeq M_\infty$, and hence
it follows from Lemma~\ref{l-5.B}  that $M$ is cell-like (and hence is acyclic), because by the above mentioned
Brown's theorem each $s_A(M)$ has this property.
Here assertion~2) of Lemma~\ref{l-5.B} applies, because the set of all finite sequences of natural numbers is well known to be countable. This completes
the proof of Lemma~\ref{l-5.2}. \hfill $\Box$

\smallskip

The answer to the converse question about monotone path-connectedness of $B$-acyc\-lic ($P$-acyc\-lic) sets in the general case is not known.
In this connection, note that for any $n\ge\nobreak 3$ one may construct an example of a~finite-dimensional space $X_n$ to contain
a~non-monotone path-\allowbreak con\-nec\-ted Chebyshev set (which is \textit{a~fortiori} a~sun and a~$P$- and $B$-acyclic set).
Indeed, consider spaces~$X_n$ with the property $\overline{\operatorname{ext}} \,S^*=S^*$.
Phelps~\cite{Phelps} has shown that  $\overline{\operatorname{ext}}
\,S^*=S^*$ holds for a~given space $X_n$ if and only if
each convex bounded closed subset of $X_n$ is representable as the intersection of closed balls
(in other words, $X_n\in (\textit{MIP})$; that is, satisfies the Mazur intersection property); as a~corollary,
in such a~space the monotone path-connectedness of a~closed set is equivalent to its convexity. Further, for any $n\ge 3$,
Tsar'kov~\cite{Tsa84} has constructed an example of a~space $X_n'$ with
 $\overline{\operatorname{ext}} \,S^*=S^*$ containing
an \textit{unbounded} nonconvex Chebyshev
set~$M'$ (any \textit{bounded\/} Chebyshev set in such $X_n'$ being convex). Thus, $M'$~serves
as an example of a~\textit{non-monotone path-connected
$B$-acyc\-lic {\rm (}$P$-acyc\-lic{\rm )} set  {\rm (}Chebyshev sun}).

In connection with this problem we point out that Brown \cite{Brown1987} has introduced an important class of normed linear spaces---the
so-called (\textit{BM})-spaces (see also \cite {Brown1988}, \cite{Brown-2002it}, \cite{FR}), which proved very natural
in the problem of  $\mcc$-connectedness of suns.
The (\textit{BM})-spaces  (in particular, $\ell^\infty(n)$ and $c_0$) are `good'  in having the property that
any boundedly compact  sun in such a~space is $\mcc$-connected (see \cite{Brown1987}, \cite{Ali05})
and hence, is monotone path-\allowbreak con\-nec\-ted.
By Vlasov's theorem, a~$P$-acyclic boundedly compact set is always a~sun, and hence in (\textit{BM})-spaces $P$-acyclicity
implies monotone path-connectedness. We also note~\cite{Brown2003} that in a~finite-dimensional \textit{polyhedral} space $X_n$
each sun is monotone path-\allowbreak con\-nec\-ted if and only if $X_n\in (\textit{BM})$.
One knows that $\ell^1(n)\notin (\textit{BM})$, $n\ge 3$, and hence in~$\ell^1(n)$ there exists a~non-monotone path-\allowbreak con\-nec\-ted sun. However, it is
unknown whether such a~sun is $P$-acyclic (or at least $P$-connected).
\medskip

The author is grateful to I.\,G.~Tsar'kov, A.\,A.~Vasil'eva, J.-D.~Hardtke, S.\,A.\,Bo\-ga\-ty\u\i\  and
U.\,Kh.~Karimov for discussions.


\begin{Biblio}
\bibitem{Brown1987} A.\,L.~Brown, ``Suns in normed linear spaces which are finite dimensional,'' Math. Ann. {\bf  279}, 1987, 87--101.

\bibitem{BDL} B.~Brosowski, F.~Deutsch, J.~Lambert, P.\,D.~Morris, ``Chebyshev sets which are not suns,'' Math.
Ann., {\bf 212} (1), 1974, 89--101.

\bibitem{Ali06}  A.\,R.~Alimov, ``Monotone path-connectedness of Chebyshev sets in the space $ C(Q)$,'' Sb. Math., {\bf 197} (9), 2006, 1259--1272.

\bibitem{Ali05}  A.\,R.~Alimov, ``Connectedness of suns in the space $ {c_0}$,'' Izv. Math., {\bf 69} (4), 2005, 651--666

\bibitem {A13FPM}   A.\,R.~Alimov, ``Local connectedness of suns in normed linear spaces,'' Fundam. Prikl. Mat., {\bf 17} (7), 2012, 3--14.

\bibitem{Vla73}  L.\,P.~Vlasov, ``Approximative properties of sets in normed linear spaces,'' Russian Math. Surveys, {\bf 28} (6), 1973, 1--66.

\bibitem{Ali12} A.\,R.~Alimov, ``A monotone path connected Chebyshev set is a sun,'' Math. Notes, {\bf 91} (2), 2012, 290--292.

\bibitem{Ali12EMJ} A.\,R.~Alimov, ``Monotone path-connectedness of $R$-weakly convex sets in spaces with linear ball embedding,'' Eurasian Math. J. {\bf 3} (2), 2012, 21--30.

\bibitem{FR}  C. Franchetti, S. Roversi, ``Suns, $M$-connected sets and $P$-acyclic sets in Banach spaces,'' Preprint no.~50139, Instituto di Matematica Applicata ``G.~Sansone'' (1988), 1--29.

\bibitem{Mel09}  S.\,A.~Melikhov, ``Steenrod homotopy'', Russian Math. Surveys, {\bf 64} (3), 2009, 469--551.

\bibitem{Massey}  W.\,S.~Massey, \textit{Homology and cohomology theory}, Marcel Dekker, New York, 1978.

\bibitem{KarRe02}  K.~Eda, U.\,H.~Karimov, D.~Repov\v s, ``On (co)homology locally connected spaces,'' Topol. Appl. {\bf 120}, 2002, 397--401.

\bibitem{Gorn}  L. G\'orniewicz, \textit{Topological fixed point theory of multivalued mappings}, Springer, 2006.

\bibitem{Gorn00}  L. G\'orniewicz, ``Topological structure of solution sets: current results,'' Arch. Math. (Brno), {\bf 36}, 2000, 343--382.

\bibitem{Dragoni}  R.~Dragoni, J.\,W.~Macki, P.~Nistri, and P.~Zecca, \textit{Solution sets of differential equations in abstract spaces}, Longman, Harlow, 1996.

\bibitem{Andres}  J.~Andres, G.~Gabor, L.~G\'orniewicz, ``Acyclicity of solution sets to functional inclusions,'' Nonlinear Analysis, {\bf 49}, 2002, 671--688.

\bibitem{HbTFT}  W. Kryszewski, ``On the existence of equilibria and fixed points of maps
under constraints,'' in: \textit{Handbook of topological fixed point theory},  Eds. R.\,F.~Brown, M.\,Furi  et al., Springer, 2005.

\bibitem {Papag}  Sh.~Hu, N.\,S.~Papageorgiou, \textit{Handbook of multivalued
analysis.} Vol. II: Applications, Kluwer, Dordrecht, 2000.

\bibitem {Brown1986}  A.\,L.\,Brown, ``Chebyshev sets and the shapes of convex bodies,'' in: \textit{Methods of Functional Analysis in Approximation
Theory}, Proc. Int. Conf., Indian Inst. Techn. Bombay,
16--20.XII.1985, Eds. C.\,A.~Micchelli, Birkh\"auser, Basel, 1986,  97--121

\bibitem{BD74}  B.~Brosowski, F.~Deutsch, ``On some geometric properties of suns,'' J. Approx. Theory, {\bf 10} (3), 1974, 245--267.

\bibitem{Giles00}   J.\,R.~Giles, ``The Mazur intersection problem,'' J. Convex Anal., {\bf 13} (3--4), 2006, 739--750.

\bibitem {Dancer}   E.\,N.~Dancer, B.~Sims, ``Weak star separability, Bull. Austral. Math. Soc.,'' {\bf 20}, 1979, 253--257.

\bibitem {Aviles}   A.~Avil\'es, G.~Plebanek, J.~Rodr\'\i guez, ``A weak$^*$ separable $C(K)^*$ space whose unit ball is not weak$^*$ separable,'' Trans. Amer. Math. Soc. (to appear)

\bibitem{Nygaard}  O.~Nygaard, ``A remark on Rainwater s theorem,'' Annales Math. Inform., {\bf 32}, 2005, 125--127.

\bibitem{Kal}    O.\,F.\,K. Kalenda, ``($I$)-envelopes of unit balls and James' characterization of reflexivity,'' Studia Math., {\bf 182}, 2007, 29--40.

\bibitem{Fabian}  M.~Fabian, P.~Habala, P.~H\'ajek, V.~Montesinos, V.~Zizler, \textit{Banach space theory. The basis for linear and nonlinear analysis}, New York, Springer, 2011.

\bibitem {TsarEC}  I.\,G.~Tsar'kov, ``Geometrical approximation theory in the papers by N.V.~Efimov and S.B.~Stechkin,'' Sovr. Probl. Mat. Mech., 2011, {\bf 4} (2), 57--78.

\bibitem{BalIva09}  M.\,V.~Balashov, G.\,E.~Ivanov, ``Weakly convex and proximally smooth sets in Banach spaces,'' Izv. Math.
{\bf 73} (3), 2009, 455--499.

\bibitem {Brown1988}   A.\,L.~Brown, ``On the connectedness properties of suns in finite         dimensional spaces,'' Proc. Cent. Math. Anal. Aust. Natl. Univ., {\bf 20}, 1988, 1--15.

\bibitem{Hu}  Sze-Tsen Hu, \textit{Theory of retracts}, Wayne State University Press, Detroit, 1965.

\bibitem{Men}  K.~Menger, ``Untersuchungen \"uber allgemeine Metrik,'' Math. Ann. {\bf 100}, 1928, 75--163.

\bibitem{Kirk}  K.~Goebel, W.\,A.~Kirk, \textit{Topics in metric fixed point theory}, Cambridge Univ.\ Press, Cambridge, 1990.

\bibitem{EGT}  K.\,P.~Hart, J.~Nagata, J.\,E.~Vaughan, \textit{Encyclopedia of general topology}, Elsevier, Amsterdam, 2004.

\bibitem{Gabor}  G.~Gabor, ``On the acyclicity of fixed point sets of multivalued maps,'' J. Juliusz Schauder Center, {\bf 14}, 1999, 327--343.

\bibitem{Dug}  J.~Dugundji, \textit{Topology}, Allyn and Bacon, Boston, MA, 1966.

\bibitem{Phelps}  R.\,R.~Phelps, ``A~representation theorem for bounded convex sets,'' Proc. Amer. Math. Soc., {\bf 11},  1960, 976--983.

\bibitem{Tsa84} I.\,G.~Tsar'kov, ``Bounded Chebyshev sets in finite-dimensional Banach spaces,'' Math. Notes {\bf 36} (1), 1984, 530--537.

\bibitem{Brown-2002it}   A.\,L.\,Brown, ``On the problem of characterising suns in finite dimensional spaces,'' Rend. Circ. Math. Palermo Ser.~II, {\bf 68}, 2002, 315--328.

\bibitem{Brown2003}    A.\,L.\,Brown , ``Suns in polyhedral spaces,'' in: \textit{Seminar of Mathem. Analysis, Proceedings},  Eds D.\,G.\,\'Alvarez, G.~Lopez Acedo and R.\,V.\,Caro,
Univ. Malaga and Seville    (Spain), Sept. 2002 -- Feb. 2003, Universidad de Sevilla, Sevilla  (2003), 139--146.
\end{Biblio}
\end{document}